\documentclass[11pt,amsart]{article}

\usepackage{epsf}

\usepackage{amssymb,amsthm,verbatim}
\font\goth=eufm10 scaled 1200

\newtheorem{prop}{Proposition}[section]
\newtheorem{thm}{Theorem}
\newtheorem{lem}[prop]{Lemma}
\newtheorem{cor}{Corollary}
\newtheorem{conj}[prop]{Conjecture}

\theoremstyle{definition}

\newtheorem{defn}[prop]{Definition}
\newtheorem{rem}[prop]{Remark}

\newtheorem{ques}[prop]{Question}

\newcommand{\bmat}{\left ( \begin{matrix} }
\newcommand{\emat}{\end{matrix} \right ) }

\newcommand{\ben}{\begin{enumerate}}

\newcommand{\een}{\end{enumerate}}

\newcommand{\blem}{\begin{lem}}

\newcommand{\elem}{\end{lem}}

\newcommand{\bcl}{\begin{clm}}

\newcommand{\ecl}{\end{clm}}

\newcommand{\bthm}{\begin{thm}}

\newcommand{\ethm}{\end{thm}}

\newcommand{\bq}{\begin{ques}}
\newcommand{\eq}{\end{ques}}

\newcommand{\bpr}{\begin{prop}}

\newcommand{\epr}{\end{prop}}

\newcommand{\bco}{\begin{cor}}

\newcommand{\eco}{\end{cor}}

\newcommand{\bcon}{\begin{conj}}

\newcommand{\econ}{\end{conj}}

\newcommand{\bde}{\begin{defn}}

\newcommand{\ede}{\end{defn}}

\newcommand{\bex}{\begin{exa}}

\newcommand{\eexa}{\end{exa}}

\newcommand{\bexe}{\begin{exe}}

\begin{document}

\title{On the Shortest Identity in Finite Simple Groups of Lie Type}

\author{Uzy Hadad}
\date{\today}                   
\maketitle                      


\begin{abstract} 
We prove that the length of the shortest identity in a finite simple group of Lie type of rank $r$ defined over $\mathbb{F}_q$, is bounded (from above and below) by explicit polynomials in $q$ and $r$.

\end{abstract}

\section{Introduction}

Let $F_\infty$ be the free group on a countable number of
generators and let $w$ be a non-trivial reduced
word in $F_\infty$ on $k$ generators, say $w=w(x_1,...,x_k)$.
Given a group $G$ we say that $w$ is an \textit{identity}, or
a \textit{law}, for $G$, if
$w(g_1,...,g_k)=1$ for every $(g_1,...,g_k)\in
G^k$. We denote by $\alpha(G)$ the length of the
shortest identity of the group $G$ (on any number of
generators).

Let $\mathbb{F}_q$ be the finite field with $q$ elements, where $q=p^m$ is a prime power, and let $G$ be a finite simple group of Lie type over $\mathbb{F}_q$ of rank $r$. Here, if $G$ is untwisted,
$r$ is the rank of the ambient simple algebraic group, while we define $r$ for the twisted groups as follows (see \cite[Sec. 13.1]{Car} and also
\cite[Prop. 2.3.2]{GLS}):

\begin{center}
  \begin{tabular}{ |l | c| c |c|c|c|c|c|c| r | }
    \hline
type   & ${}^2A_{2d}$ & ${}^2A_{2d-1}$  & ${}^2D_{d+1}$ & ${}^3D_4$ & ${}^2E_6$ & ${}^2F_4$   & ${}^2G_2$  & ${}^2B_2$   \\ \hline
    rank  & $d$ & $d$ & $d$ & $2$ & $4$ & $2$& $1$& $1$ \\ \hline

    \hline
  \end{tabular}
\end{center}

Let $q^*(G)$ be the number of elements of the field where the group $G$ is realized. For example, $q^*(A_n(q))=q$, $q^*({}^3D_4(q))=q^3$ and $q^*({}^2A_{2d}(q))=q^2$. See Section \ref{not_def} for more details.

The main result of this paper is the following:

 \bthm
\label{thm:main} Let $G$ be a finite simple group of Lie
type over $\mathbb{F}_q$ of rank $r$. Then the length of the shortest identity of
$G$ satisfies
 $$\frac{q^*(G)^{\frac{r}{4}}-1}{3} \leq \alpha(G) < (31 r+2)^3q^{31 r}.$$
Furthermore, we give an explicit construction
of an identity of length less than the upper bound.

 \ethm

In the case of $G=A_d(q)=PSL_{d+1}(q)$ we give a more
precise statement, and then deduce Theorem \ref{thm:main} from
this special case. Indeed, we prove the following:

 \bthm \label{thm:sl}
The length of the shortest identity of $A_d(q)$ satisfies
 $$ \frac{q^{\lfloor\frac{d+1}{2}\rfloor}-1}{3} \leq \alpha(A_d(q)) <  (d+3)^2dq^{d+1}.$$
Furthermore, we give an explicit construction
of an identity of length less than the upper bound.
\ethm

In particular, Theorem \ref{thm:sl} improves a result of Gamburd et. al. \cite[Prop. 11]{GHSSV}
 which states that $\alpha(SL_2(p)) \geq c\frac{p}{\log p}$
for some constant $c$.

The proof of Theorem \ref{thm:main} is
based on the fact that if $G$ is a finite simple group of Lie type over $\mathbb{F}_q$ (with the exception of Suzuki group) then
there exists two positive integers $c_1$ and $c_2$, which depend only on the type and the rank of $G$,
such that $$A_1(q^{c_1}) \leq G \leq A_{c_2}(q).
$$

Now since an identity of a group is
inherited by its subgroups, we observe that Theorem \ref{thm:main} follows from Theorem \ref{thm:sl}.

\begin{rem}
If a group $G$ is a central extension of a group $Z$ by $H$( i.e
$H=G/Z$), then if a word $1\neq w(x_1,...,x_k) \in F_\infty$ is an
identity of $H$ then it is easy to check that the word
$[w,x_{k+1}] \in F_\infty$ is also an identity of $G$. Furthermore,  if $1 \neq
w(x_1,...,x_k) \in F_\infty$ is an identity of $G$ then it is an
identity of $H$. Therefore we obtain similar results for quasisimple groups.
\end{rem}

Let $F_k$ be the free group on $k$ generators $x_1,...,x_k$. Given a word  $w(x_1,...,x_k) \in F_k$ we define $l(w)$
 to be the length of $w$.
 Let $1 \neq w \in F_k$ be an identity for a finite group $G$
and let $S=\{g_1,...,g_k\}\subset G$. Then starting at any vertex in the Cayley graph $Cay(G,S)$, the walk $w(g_1,...,g_k)$
is a closed walk of length $l(w)$. Hence in some
sense, $\alpha(G)$ is a universal girth of the Cayley
graphs of the group.

Recently there has been great interest in the study of word
maps in finite simple groups and in algebraic groups (see
\cite{Lar,LS,NS1,Sh}). Let $w=w(x_1,...,x_k)$ be a non-trivial
word in $F_k$, where $k \geq 1$, and let $G$ be a group. We define
$$w(G)=\langle w(g_1,...,g_k):(g_1,...,g_k) \in G^k\rangle.$$

The next result follows immediately from Theorem \ref{thm:main} (since the hypothesis on $l(w)$ implies that
$w(G)$ is a non-trivial normal subgroup of $G$).
 \bco Let $G$ be a finite simple group of Lie type over $\mathbb{F}_q$ of rank $r$, and let $w \in F_k$ be a non-trivial word with  $l(w) < \frac{q^*(G)^{\frac{r}{4}}-1}{3}$, where $k \geq 1$. Then $w(G)=G.$

\eco

\subsection{Notation and Definitions} \label{not_def}

 Let $G$ be group, for every
$g,h, \in G$ we define $g^h=hgh^{-1}$ and
$[g,h]=ghg^{-1}h^{-1}$.
 In case $G$ is finite, we write $exp(G)$ for the exponent of $G$.
In addition, for an element $g\in G$ we define $ord(g)$ to be
the order of $g$. For $x \in \mathbb{R}$ we write $\lfloor x
\rfloor$ for the greatest integer less than or equal to
 $x$.

The untwisted groups $A_d(q),B_d(q),C_d(q),D_d(q),E_6(q),E_7(q),E_8(q),F_4(q),G_2(q)$ are
realized over the finite field $\mathbb{F}_q$.
The twisted groups
$${}^2A_{2d}(q),{}^2A_{2d-1}(q),{}^2D_{d+1}(q),{}^3D_4(q),{}^2E_6(q),$$
are realized over finite fields with
$q^2,q^2,q^2,q^3,q^2$
elements respectively, and the Ree and Suzuki groups ${}^2F_4(q),{}^2G_2(q),{}^2B_2(q),$
are realized over finite fields with
$q=2^{2n+1},q=3^{2n+1},q=2^{2n+1}$ elements respectively
(see \cite[ch. 13.] {Car}).

 We shall not go into further details about the
structure and construction of finite simple groups of Lie type, we
refer the reader to the book of Carter \cite{Car} which is our
main reference (see also \cite{GLS}).

\section{Previous Work} \label{pre_work}
A \textit{variety} of groups is a class of groups that satisfy
a given set of laws (see \cite[ch. 1]{Ne}).  By Birkhoff's theorem \cite{Be} each variety $\hbox{\goth  B}\hskip 2.5pt$
is defined by a suitable set of words $B \subseteq F_\infty$, that is, $\hbox{\goth  B}\hskip 2.5pt$ consists of all the groups
$G$ on which $w(g_1,...,g_k)=1$ holds for each word $w(x_1,...,x_k)\in B$ and for each set of elements $g_1,...,g_k \in G$.

For example the variety which is
defined by the law $$\{[x_1,x_2]=x_1x_2x_1^{-1}x_2^{-1}=1\},$$ is the abelian groups. In the language of varieties, the main aim
of this paper is to study the length of the shortest
law in the variety which is generated by all the laws in a finite
simple group of Lie type.

In the 1960's, various papers were written on varieties of groups,
which were mainly concerned with their qualitative properties. The
most notable contribution is Hanna Neumann's book \cite{Ne}. In this
book, she raised (\cite[p. 166]{Ne}) the following question:

 \begin{quote} \textit{Is there a law which is satisfied in an infinite number of
non-isomorphic non-abelian finite simple groups?}
\end{quote}

G. A. Jones \cite{Jon} gave a negative answer to this question,
but his proof does not give an explicit
bound on the length of the shortest identity in each family of finite simple groups of Lie type (except for the case of Suzuki groups).
Our work can be thought of as a
quantitative version of Jones's results.

A \textit{basis} for a variety $\hbox{\goth  B}\hskip 2.5pt$ is a set of laws such that its closure is $\hbox{\goth  B}\hskip 2.5pt$ (for the definition  of `closure', see \cite[ch. 1]{Ne}). Oates and Powell in \cite{OP} proved that every variety
generated by a finite group has a finite basis.
It is clear that if we know a basis for the variety generated by the laws in a given group $G$,
then the minimum length of a law in this basis, is an upper
bound for the shortest identity of the group.

In the literature a few attempts have been made to find an explicit basis for the
 variety generated by a finite non-abelian simple group. J.
Cossey and S. Macdonald \cite{CM} gave a finite basis for the set
of laws in $PSL_2(5)$,  and this was extended in \cite{CMS} to
$PSL_2(p^n)$ with $p^n\leq 11$. B. Southcott \cite{So1,So2} gave a
basis for the family $PSL_2(2^n)$, but the length of each element
in his basis is greater than the upper bound of the shortest
identity we state in Theorem \ref{thm:sl}.

Let $\mathcal{X}$ be any infinite set of groups. A group $G$
is said to be \textit{residually} $\mathcal{X}$, if for every $1\neq g \in
G$ there exists an epimorphism $\varphi$ from $G$ to some
$H \in \mathcal{X}$ such that $\varphi(g)\neq 1$.

Suppose that the free group on $k$ generators $F_k$ is residually $\mathcal{X}$, and for a given group $H \in \mathcal{X}$
we can determine the maximal length of a non-trivial word $w$ in $F_k$ such that there exists an  epimorphism $\varphi$ from $F_k$
to $H$ with $\varphi(w)\neq 1$. Then this gives a lower bound on the length of the shortest identity  (on $k$ generators) in $H$.

 W. Magnus \cite[p. 309]{Ma} raised the following related problem:

\begin{quote}\textit{Let $\mathcal{X}$ be any infinite set of non-abelian finite
simple groups. Is the free group $F_k$ on $k\geq 1$ generators
residually $\mathcal{X}$?}
\end{quote}

\noindent T. Weigel \cite{We1,We2,We3} gave a complete answer to this
question. From his proof, we can conclude that in the case of a
classical simple group $G$ over $\mathbb{F}_q$  where $q=p^m$ is a prime
power, the length of the shortest identity (on two generators) is
at least $p$, and for the exceptional
groups $E_6(q),E_7(q),E_8(q),G_2(q),F_4(q)$ and the twisted groups
${}^3D_4(q),{}^2E_6(q),{}^2F_4(q)$, the length of the shortest identity (on two generators)
is at least $\log q$.
Weigel's work involves `identities' on two generators
(and their inverses) of the group; in this paper we consider
identities on any number of elements of the group (not necessarily
generators).

 \section {Proof of Theorem \ref{thm:sl}} \label{proof_thm_sl}

\subsection{The Upper Bound}

In this section we give an explicit construction of an identity for the group $SL_n(q)$ and this gives the upper bound stated in Theorem \ref{thm:sl}.
This construction is based on the exponent of the group.

\blem \label{lem_exp} Let $q=p^f$ where $p$ is a prime. Then the
exponent of $SL_n(q)$ is
$$p^e\cdot lcm[q-1,q^2-1,...,q^{n-1}-1,\frac{q^n-1}{q-1}],$$
where $e$ is the minimal positive integer such that $p^e \geq n$.

 \elem

\begin{proof}

Let $x\in GL_n(q)$ be a non-trivial element. Write $x=x_s\cdot
x_u$ in Jordan form, where $x_u\in GL_n(q)$ is unipotent, $x_s \in
GL_n(q)$ is diagonalizable over a splitting field $\mathbb{F}_{q^r} (r
\leq n)$ and $[x_s,x_u]=1$. Here $ord(x_u)$ is a power of $p$, while $ord(x_s)$ is coprime to $p$.

Suppose $x_u=1$. Then $x$ is diagonalizable
over $\mathbb{F}_{q^r}$, and therefore
$x^{q^r-1}=x_s^{q^r-1}=1$. Furthermore, for every $1 \leq j\leq
n$, there exists an element in $GL_n(q)$ of order $q^j-1$. This
follows from the fact that the field $\mathbb{F}_{q^j}$ can be
considered as a vector space over $\mathbb{F}_q$ of dimension $j$.
Take a generator of the cyclic group $\mathbb{F}_{q^j}^*$; it is of order
$q^j-1$ and it acts as a linear transformation on
$\mathbb{F}_q^j$, hence it belongs to $GL_n(q)$.

Now suppose $x_u$ is non-trivial and let $m=ord(x)$. Then $x_s^mx_u^m=1$ since $[x_s,x_u]=1$,
so $m$ is divisible by $ord(x_s)$ and $ord(x_u)$.


Now every unipotent matrix $x_u$ can we written as $x_u=1+N$,
where $1$ is the identity matrix in $GL_n(q)$ and $N$ is an
$n\times n$ nilpotent matrix over $\mathbb{F}_{q}$. For every $i
\in \mathbb{N}$, we have $x_u^{p^i}=(1+N)^{p^i}=1+N^{p^i}$.  Let
$e$ be the minimal positive integer such that $p^e \geq n$. Then
it is easy to check that $p^e$ is the exponent of the upper
unitriangular matrices in $GL_n(q)$, and there exists a unipotent matrix in $GL_n(q)$ with this order.

For $SL_n(q)$ all the above holds  except for the fact that the
maximal order of an element in $SL_n(q)$ is $\frac{q^n-1}{q-1}$
and not $q^n-1$ as in $GL_n(q)$.
\end{proof}

\begin{rem}
There exists a constant $c$ such that $exp(SL_n(q))\geq
q^{cn^2}$ (see \cite[Lemma 2.3]{BMP} and the discussion afterward).
Thus our bound in Theorem \ref{thm:sl} is shorter than the length
of the exponent identity $x^{exp(SL_n(q))}$.
\end{rem}

 \blem \label{lem_cons_word} Let $G$ be a group and let $w_1,...,w_m$ be
distinct non-trivial power-words in one variable in the free group
$F_{k+1}$ on $k+1$ generators $x_1,...,x_{k+1}$ where
$k=2^{\lfloor \log_2m \rfloor}$. Suppose that $l(w_1) \geq l(w_2)
\geq ... \geq l(w_m)$ and that for each element $g \in
G$, there exists some $i$ such that $w_i(g)=1$. Then
there exists a non-trivial word $w\in F_{k+1}$ of length at most
$4m^2(l(w_1)+1)$ which is an identity in $G$. \elem

\begin{proof}

For $1 \leq i \leq \lfloor\frac{m}{2}\rfloor$, set
$u_i=[w_{2i-1}(x_1),w_{2i}(x_1^{x_{i+1}})]$.
 If $m$ is odd then set
$$u_{\lfloor\frac{m}{2}\rfloor+1}=[w_m(x_1),x_{\lfloor\frac{m}{2}\rfloor+2}],$$
else if $m$ is even set
$$u_{\lfloor\frac{m}{2}\rfloor+1}=[x_1,x_{\lfloor\frac{m}{2}\rfloor+2}].$$
For $\lfloor\frac{m}{2}\rfloor+2 \leq i \leq k$, set
$u_i=[x_1,x_{i+1}]$. Let $j=2^e$ for some positive integer $e \leq
\lfloor \log_2m \rfloor$, and define a recursive function $f$ in the
following way:
$$f(u_1,...,u_{2^{e-1}},...,u_j)=[f(u_1,...,u_{2^{e-1}}),f(u_{2^{e-1}+1},...,u_j)]$$
and $f(u_i,u_{i+1})=[u_i,u_{i+1}]$.

Let $w=f(u_1,...,u_k).$ It is clear that $w$ is a non-trivial word
in $F_{k+1}$ (since we have a new letter $x_{j+1}$ in each $u_j$,
which appears only in $u_j$). To show that $w$ is an identity of
$G$ it is enough to note that $x$ and $x^y$ have the same
order, hence at least one of the commutators in the expression for
$w$ collapses and hence so does the whole word.

 Now since $w_i(x_1^{x_j})=x_jw_i(x_1)x_j^{-1}$, the length of the word $[u_1,u_2]$
is at most $2^4\cdot l(w_1)+2^4$. By induction on $e$, since $l(w_1)\geq l(w_i)$ for all $i$, we get that
the length of the word $f(u_1,...,u_{2^e})$  is at most
$$2^{2(e+1)} l(w_1)+ 2^{2(e+1)}.$$
Therefore, since $e \leq \lfloor \log_2m
\rfloor$, the length of $w$ is at most $$2^{2({\lfloor \log_2m
\rfloor}+1) }l(w_1)+2^{2({\lfloor \log_2m \rfloor}+1 )} \leq
4m^2l(w_1)+4m^2 =4m^2(l(w_1)+1).$$

\end{proof}

 \bpr \label{pr:sl:iden} Let $q=p^f$ be a prime power and let $G=SL_n(q)$ where $n\geq2$. Then there exists an
identity in $G$ of length at most $ (n+2)^2 p^e q^{n-1}$,
where $e$ is the minimal positive integer such that $p^e \geq
n$.\epr

\begin{proof}
Following the proof of  Lemma \ref{lem_exp}, every element $g
\in G$ satisfies at least one of the following words:
$$X^{p^e\frac{q^n-1}{q-1}},X^{p^e(q^{n-1}-1)},X^{p^e(q^{n-2}-1)},...,X^{p^e(q-1)}.$$
If $g$ satisfies the first word, then it also satisfies
$X^{\frac{q^n-1}{q-1}}$ since it has distinct eigenvalues in the
splitting field $\mathbb{F}_{q^n}$, hence it is diagonalizable over $\mathbb{F}_{q^n}$. Now
for every $ i \in \mathbb{N}$, $q^i-1$ divides $q^{2i}-1$,
therefore every $g \in G$ satisfies at least one of the
following words (ordered in decreasing length):
$$w_1=X^{p^e(q^{n-1}-1)},w_2=X^{\frac{q^n-1}{(q-1)}},w_3=X^{p^e(q^{n-2}-1)},...,w_{\lceil\frac{n+1}{2}\rceil}=X^{p^e
(q^{\lceil\frac{n+1}{2}\rceil}-1)}.$$

\noindent Now set $m=\frac{n+2}{2}$. The result follows by Lemma \ref{lem_cons_word}.
 \end{proof}

If $w$ is an identity in the group $G$, it is easy to see
that $w$ is also an identity for every quotient of $G$. Hence
the identity we construct in Proposition \ref{pr:sl:iden} holds
for $A_{n-1}(q)=PSL_n(q)$, and this gives the upper bound in
Theorem \ref{thm:sl}.

\subsection {The Lower Bound}

 \blem \label{len_even} Let $G$ be a finite group and $w$ an identity of
 $G$. If $l(w) < exp(G)$, then $l(w)$ is even.

 \elem
 \begin{proof}
Let $w=w(x_1,...,x_k)$ be an identity in $G$ on $k$ elements.
If $l(w)$ is odd, then there exists $ 1 \leq i \leq
k$ such that the sum of the exponents of $x_i$ appearing in $w$
is odd. Hence if we set $x_j=1$ for all $j \neq i$, we get a power-word of
length less than $exp(G)$, a contradiction.

 \end{proof}

We start with the case $A_1(q)=PSL_2(q)$.
 \blem \label{lem:sl_2}
Let $G=PSL_2(q)$, then $\alpha(G) \geq \frac{q-1}{3}$.

 \elem
\begin {proof}

If $q$ is even then $PSL_2(q)=SL_2(q)$. If $q$ is odd then the
group $PSL_2(q)$ is isomorphic to $SL_2(q)/\mathbb{Z}_2$. In any
case, if $w \in F_k$ is an identity of $PSL_2(q)$, then it is
immediate to check that $w^2$ is an identity of $SL_2(q)$. So it
is enough to prove that the length of the shortest identity of
$SL_2(q)$ is at least $\frac{2}{3}(q-1)$.

  Let $$u(t)=\left(\begin{array}{cc}
  1 & t  \\
  0 & 1\end{array} \right),
  \tau= \left(\begin{array}{cc}
  0 & -1  \\
  1 & 0\end{array} \right)
 ,h(\lambda) =  \left(\begin{array}{cc}
  \lambda & 0  \\
  0 & \lambda^{-1} \end{array} \right)$$
where $t \in \mathbb{F}_q$ and $\lambda \in \mathbb{F}_q^*$. From the
Bruhat decomposition for $SL_2(q)$, every element $g \in SL_2(q)$
has a unique expression in one of the following forms (see
\cite[Cor. 8.4.4]{Car}):

$$g=u(a)h(\lambda) \mbox{   or   } g=u(b)h(\gamma)\tau u(c).$$

Suppose $w=w(x_1,...,x_k)$ is a non-trivial reduced word in $F_k$
of length $l$ which is an identity of $SL_2(q)$. Since we are interested in deriving a
lower bound, we may assume that $l$ is less than
$exp(SL_2(q))$ (see Lemma \ref{lem_exp}), so $l$ is even
by Lemma \ref{len_even}.

Let
$$M=\{a_i,b_i,c_i,\lambda_i,\gamma_i:1 \leq i \leq k \}$$ be a set of
independent commuting indeterminates over $\mathbb{F}_q$, and let
$$X_i=u(a_i)h(\lambda_i) \mbox{     and     } Y_i=u(b_i)h(\gamma_i)\tau
u(c_i).$$ Then $X_i$ and $Y_i$ are matrices with entries in
$\mathbb{F}_q[M]$ and it is immediate to verify that the entries
of
$$\lambda_iX_i,\gamma_iY_i,\lambda_iX_i^{-1} \mbox{  and  } \gamma_iY_i^{-1}$$ are
polynomials of degree at most $2$ in the variables in $M$.

Let $n_i$ be the sum of the moduli of the exponents of $x_i$
appearing in $w$ (for example if $w=x_1x_2x_1^{-1}x_2^{-1}$ then
$n_1=n_2=2$) and let $I_2$ be the identity matrix in $SL_2(q)$. For
any $Z_i \in \{X_i,Y_i\}$, the matrix

$$C(Z_1,...,Z_k)=\prod \limits_{i=1}^k\beta(Z_i)^{n_i}(w(Z_1,...,Z_k)-I_2)$$

\noindent where $\beta_i(X_i)= \lambda_i$ and $\beta(Y_i) = \gamma_i$,
 has entries in $\mathbb{F}_q[M]$ having degree at most
$2n_i$ in each of the variables $a_i,b_i,c_i,\lambda_i,\gamma_i$.

\noindent Now there are two cases to consider:
\begin{enumerate}
\item[(i)]  For all the substitutions of $Z_i$ by $X_i$ or $Y_i$ we get $C(Z_1,...,Z_k)=0$.
 \item[(ii)]  There is a substitution
$(Z_1,...,Z_k)$ such that $C(Z_1,...,Z_k)\neq 0$.
\end{enumerate}

First let us consider case (i). Let $K$ be the algebraic closure
of $\mathbb{F}_q$. Since for every substitution $C(Z_1,...,Z_k)$ is zero, we deduce that $w$ is a law on $SL_2(K)$. For
every $n\in \mathbb{N}$, $SL_2(K)$ has a subgroup isomorphic to
$SL_2(q^n)$, hence we obtain a law for the infinite family
$\{SL_2(q^n)\}_{n=1}^\infty$. But the main theorem of Jones
\cite{Jon} states that there is no law which is satisfied by an
infinite family of non-abelian finite simple groups, a
contradiction.

Now let us consider case (ii). Let $(Z_1,...,Z_k)$ be a substitution such that the word
$w(Z_1,...,Z_k)-I_2$ is not zero.
 Decompose $w(Z_1,...,Z_k)$ as a product of two words:
$$w(Z_1,...,Z_k)=w_1(Z_1,...,Z_k)w_2(Z_1,...,Z_k)$$ where
$w_1(Z_1,...,Z_k)$ is a word of length $\frac{1}{2}l$ (recall that $l$ is
even). Let $T$ be the matrix
$$T= \prod \limits_{i=1}^k\beta(Z_i)^{n_i}(w_1(Z_1,...,Z_k)-w_2(Z_1,...,Z_k)^{-1}).$$

\noindent Then there is an entry $T_{i,j}$  for some $1 \leq i,j
\leq 2$, which is not formally zero.  Let $M_1 \subseteq M$ be the
set of indeterminates appearing in $T_{i,j}$. Then  $T_{i,j}$ is a
polynomial in $|M_1|$ variables of degree at most $\frac{3l}{2}$.
Recall that if $f$ is a polynomial in $\mathbb{F}_q[M_1]$ with $deg(f)=d\geq 0$,
then the equation $f=0$ has at most $dq^{|M_1|-1}$ solutions in
$\mathbb{F}_q^{|M_1|}$ (see \cite[Thm. 6.13]{LN}). Therefore $T_{i,j}$
has at most $\frac{3}{2}l \cdot q^{|M_1|-1}$ solutions, hence $l
\geq \frac{2}{3}(q-1)$ (the $q-1$ factor is because
$\lambda_i,\gamma_i \in \mathbb{F}_q^*$).

\end{proof}

It is easy to see that $A_1(q^n)=PSL_2({q^n})$ is a subgroup of
$A_{2n-1}(q)=PSL_{2n}(q)$ (this follows from the fact that a $2$ dimensional vector space over $\mathbb{F}_{q^n}$, can be considered as a $2n$ dimensional vector space over $\mathbb{F}_q$), thus we obtain the following
result:
 \bpr Let $G=A_d(q)$ where $d\geq 1$. Then $\alpha(G)>\frac{q^{\lfloor\frac{d+1}{2}\rfloor}-1}{3}$. \epr

This completes the proof of Theorem \ref{thm:sl}.

 \section {Proof of Theorem \ref{thm:main}} \label{proof_thm_main}

\subsection{The Lower Bound} \label{sec_lw}

In this section we will show that if $G$ is a finite simple group of Lie type of rank $r$, and $G$ is not a Suzuki group (i.e. type ${}^2B_2$), then $G$ contains $A_{r'}(q')$,
where $q'\leq q^*(G)$ and $r' \leq r$.
Therefore, using Theorem \ref{thm:sl}, we obtain a lower bound for the shortest identity of $G$. For
the Suzuki  groups we use a lemma of Jones.

\subsubsection{Untwisted groups}

 Let $G\neq A_d(q)$  be a simple untwisted group of Lie type over $\mathbb{F}_q$.
 By considering the associated Dynkin diagram, it is clear that if $G$ has rank $d$ then $G$ contains $A_{d-1}(q)$, so
 Theorem \ref{thm:sl} implies that $\alpha(G) \geq \frac{q^{\lfloor
\frac{d}{2}\rfloor}-1}{3}$.

\subsubsection{Twisted groups}

\textbf{Ree Groups:} It is known that $^2F_4(2^{2n+1})$
and $^2G_2(3^{2n+1})$ contains $A_1(3^{2n+1})$ and
$A_1(3^{2n+1})$ respectively ~\cite{Ti,Lh}.

\noindent \textbf{Suzuki groups:} The following lemma is a
quantitative version of a result of Jones.
 \blem The length of the shortest identity of $^2B_2(2^{2n+1})$ satisfies $$\alpha(^2B_2(2^{2n+1})) \geq
\frac{2^{2n}-1}{1+2^n}.$$\elem

For the proof we refer the reader to \cite[Lemma 5]{Jon}).
Although this quantitative result is not stated there explicitly,
it follows immediately from the proof (in fact, one can improve this bound
by decomposing the identity
word in the same way we did in the proof of Lemma \ref{lem:sl_2}).\\

\noindent \textbf{The twisted groups
$^2A_d(q)$,$^2D_d(q)$,$^2E_6(q)$,$^3D_4(q)$:}

The Dynkin diagrams of type $A_d,D_d,E_6$ and $D_4$ admit symmetries of order $2,2,2$ and $3$, respectively.
 The twisted groups  $^2A_d(q)$,$^2D_d(q)$,$^2E_6(q)$ and $^3D_4(q)$ are subgroups of
the untwisted groups $A_d(q^2),D_d(q^2),E_6(q^2),D_4(q^3)$, and
they are fixed by an automorphism $\sigma$ (of $A_d(q^2)$ etc.) which maps each root
element $X_\alpha(t)$ to $X_{\alpha'}(t^q)$, where $\alpha
\mapsto \alpha'$ is a symmetry of the root system and $t
\mapsto t^q$ is a field automorphism. The automorphism $\sigma$
has order $2,2,2,3$ respectively.

By inspecting the Dynkin diagrams of
type $A_{2d-1},A_{2d},D_d,E_6$ and $D_4$, together with the corresponding
automorphisms and roots relations, one can show that each of the following groups $G$ has a subgroup $H$ as given
in the following table (see also \cite[Sec. 2]{Ni}):

\begin{center}
  \begin{tabular}{ |l | c| c|c |c|c|c| r | }
    \hline
     $G$ & $^2A_2(q)$ & $^2A_{2d-1}(q)$& $^2A_{2d}(q),d >1$ & $^2E_6(q)$ & $^3D_4(q)$ \\ \hline
     $H$ & $A_1(q)$ & $A_{d-1}(q^2)$ & $A_{d-1}(q^2)$ &  $A_{2}(q^2)$ & $A_1(q^3)$\\ \hline

    \hline
  \end{tabular}
\end{center}

We give the details in the case of ${}^2A_{2d-1}(q)$, a similar argument applies in each of the remaining cases.
Let $\prod=\{\omega_1,...,\omega_{2d-1}\}$ be a base for the root system of type $A_{2d-1}$.
 Now for every $t \in \mathbb{F}_{q^2}$ and for every $ 1 \leq i <d$
the element $X_{\omega_i}(t)X_{\omega_{2d-i}}(t^q)$ is fixed
by the automorphism $\sigma$ and for every $t\in \mathbb{F}_q$ the element $X_{\omega_d}(t)$ is fixed by $\sigma$.
It is well-known fact that if $|i-j|>1$ then the root subgroups $X_{\omega_i}$ and $X_{\omega_j}$ commute, thus we get that
the subgroup $$H=\langle X_{\omega_i}(t)X_{\omega_{2d-i}}(t^q): 1 \leq i <d,t\in \mathbb{F}_{q^2}\rangle$$
 is isomorphic to $A_{d-1}(q^2)$.

So for a twisted group $G$ in the above table with rank
$r$, we deduce from Theorem \ref{thm:sl} that the shortest identity of $G$ has length at least
$\frac{(q^*(G))^{\frac{r}{4}}-1}{3}$.

The only case left is $^2D_{d+1}(q)$ (recall that the rank of
$^2D_{d+1}(q)$ is $d$). Let $n=2^k$ for some integer $k\geq 2$ and fix $d$ such that $ n \leq d+1 <2n =2^{k+1}$. Then we have
$^2D_{d+1}(q) \geq {^2D_{\frac{n}{2}}}(q^2) \geq A_1(q^n)$ (see
\cite[Table 3.5.F]{KL}). Now since $2(\frac{d}{4}) < n$,
we deduce that the shortest identity of $^2D_{d+1}(q)$ has length at least
$ \frac{q^n-1}{3} > \frac{(q^2)^\frac{d}{4}-1}{3}$.

\subsection{The Upper Bound} \label{sec_up}
In this final section we show that there exists a constant $c$ such that
every finite simple group of Lie type over $\mathbb{F}_q$ of rank $r$ is isomorphic to
a subgroup of $PSL_{cr}(q)$ or $SL_{cr}(q)$, and so the desired upper bound follows from Proposition \ref{pr:sl:iden}.

Each finite simple group of Lie type can be constructed (see \cite[ch. 4]{Car}) as a subgroup
of the automorphism group of the corresponding Lie
algebra.
In the following table we give the dimension of the Lie algebra $\hbox{\goth  g}\hskip 2.5pt$ in
terms of the untwisted Lie rank:

\begin{center}
  \begin{tabular}{ |l | c| c |c|c|c|c|c|c|c| r | }
    \hline
    \hbox{\goth  g}\hskip 2.5pt   & $A_d$ & $B_d$  & $C_d$ & $D_d$ & $E_6$ & $E_7$   & $E_8$  & $F_4$  & $G_2$ \\ \hline
    dim \hbox{\goth  g}\hskip 2.5pt  & $d(d+2)$ & $d(2d+1)$ & $d(2d+1)$ & $d(2d-1)$ & $78$ & $133$& $248$& $52$ & $14$\\ \hline

    \hline
  \end{tabular}
\end{center}

It is well known (see \cite[Ch. 11]{Car}) that each of the simple
classical groups
$$A_d(q)\mbox{,}^2A_d(q),B_d(q),C_d(q),D_d(q)\mbox{,}^2D_d(q)$$ has a matrix
representation as a subgroup of $PSL_{2d+1}(q^2)$.

\noindent In addition, the above table indicates that the exceptional groups
$$E_6(q),E_7(q),E_8(q),F_4(q),G_2(q)$$
are subgroups of $SL_{248}(q)$ (via the adjoint representation on the corresponding Lie algebra),
so for the untwisted groups we can
take $c=31$.

Each twisted group $G$ is the set of fixed points of some
automorphism of the corresponding untwisted group, hence a
subgroup of some $SL_{cr}(q)$ (or $PSL_{cr}(q))$ for some $c$. For example,
$$^2F_4(q) < F_4(q^2)<SL_{52}(q^2)<SL_{104}(q).$$
One can check from the table above that $c = 31$, and this finishes
the proof of Theorem \ref{thm:main}.

\section {Acknowledgments} This paper is part of the author's PhD thesis.
The author is grateful to his advisor Prof. Alex Lubotzky for
introducing him  to the problem and for useful discussions. Thanks to
Prof. Inna (Korchagina) Capdeboscq, Prof. Nati Linial, Prof.
Avinoam Mann, Dr. Nikolay Nikolov and Prof. Aner Shalev for useful remarks. Thanks are
also due to M. Berman for his careful reading and for the referee for
his constructive suggestions and criticisms which helped me improve this paper.


\noindent Einstein Institute of Mathematics, The Hebrew University of Jerusalem. \\
\textit{current:} Department of Mathematics, Weizmann Institute of Science. \\ 
and\\
Computer Science Division, The Open University of Israel.\\
\textit{E-mail address:} uzy.hadad@gmail.com

 \end{document}